\documentclass[letterpaper,10pt,twocolumn]{IEEEtran}

\IEEEoverridecommandlockouts        
\renewcommand{\baselinestretch}{0.99}

\title{\LARGE \bf Sparsity-Promoting Optimal Wide-Area Control of Power Networks}

\author{Florian D\"orfler, Mihailo R. Jovanovi\'c, Michael Chertkov, and Francesco Bullo%
  \thanks{This material is based in part upon work supported by the NSF grants IIS-0904501, CPS-1135819, and CMMI-09-27720. A preliminary and abbreviated version of this document has been presented in~\cite{FD-MJ-MC-FB:12s}.}%
  \thanks{Florian D\"orfler is with the Electrical Engineering Department, University of California Los Angeles, Los Angeles, CA 90095. Mihailo R.\ Jovanovi\'c is with the Department of Electrical and Computer Engineering, University of Minnesota, Minneapolis, MN 55455.  Michael Chertkov is with Theory Division \& Center for Nonlinear Studies at LANL and with New Mexico Consortium, Los Alamos, NM 87544.  Francesco Bullo is with the Department of Mechanical Engineering, University of California  Santa Barbara, Santa Barbara, CA 93106.  Email: {\tt dorfler@seas.ucla.edu}, {\tt mihailo@umn.edu}, {\tt chertkov@lanl.gov}, {\tt bullo@engineering.ucsb.edu}.}}

\usepackage{graphicx,color}
\usepackage{amsmath}
\usepackage{amssymb}
\usepackage{mathrsfs}
\usepackage{dsfont}
\usepackage[vlined,ruled]{algorithm2e}
\usepackage{subfigure}
\usepackage[noend]{algorithmic}
\usepackage{cite}

\newcommand{\subscr}[2]{{#1}_{\textup{#2}}}
\newcommand{\until}[1]{\{1,\dots,#1\}}

\newcommand{\real}{\mathbb{R}}

\DeclareMathOperator*{\minimize}{minimize}
\newcommand{\mc}{\mathcal}

\newcommand\oprocendsymbol{\hbox{$\square$}}
\newcommand\oprocend{\relax\ifmmode\else\unskip\hfill\fi\oprocendsymbol}

\renewcommand{\theenumi}{(\roman{enumi})}

\begin{document}
\maketitle
\thispagestyle{empty}
\pagestyle{empty}
 


\begin{abstract}
Inter-area oscillations in bulk power systems are {typically} poorly controllable by means of local decentralized control.  {Recent} research efforts {have been aimed at developing} wide-area control strategies {that involve} communication of remote signals. In conventional wide-area control, the control structure is fixed {\em a priori\/} typically based on modal criteria. In contrast, here we employ the recently-introduced paradigm of sparsity-promoting optimal control to simultaneously identify the optimal control structure and optimize the closed-loop performance.
To induce a sparse control architecture, we regularize the standard quadratic performance index with an $\ell_{1}$-penalty on the feedback matrix. The quadratic objective functions are inspired by the {classic} slow coherency theory and are aimed at imitating homogeneous networks without inter-area oscillations. 
We use {the New England power grid model} to demonstrate that the proposed combination of the sparsity-promoting control design with the slow coherency objectives performs almost as well as the optimal centralized control while only making use of a single wide-area communication link. 
{In addition to} this nominal performance, we also demonstrate that our control strategy yields favorable robustness margins and {that it} can be used to identify a sparse control architecture for~control~design via {alternative means}.
\end{abstract} 



\section{Introduction}
\label{Section: Introduction}

Bulk power systems typically exhibit multiple electromechanical oscillations. {\em Local oscillations\/} are characterized by single generators swinging relative to the rest of the network, whereas {\em inter-area oscillations} are associated with the dynamics of power transfers and involve groups of generators oscillating relative to each other. 
With steadily growing demand, deployment of renewables in remote areas, and deregulation of energy markets, long-distance power transfers outpace the addition of new transmission facilities. These developments lead to a maximum use of the existing network, result in smaller stability margins, and cause inter-area modes to be ever more lightly damped.
%
In a heavily stressed grid, poorly damped inter-area modes can even become unstable, 
For example, the blackout of August 10, 1996, resulted from an instability of the $0.25\textup{Hz}$ mode in the Western interconnected system~\cite{VV-LV:04}.

Local oscillations are typically damped by generator excitation control via power system stabilizers (PSSs) \cite{PK:94}. However, these {\em decentralized} control actions can interact in an adverse way and destabilize the overall system \cite{RG-KI-SL-JP-RC:93}.~Sometimes inter-area modes cannot be {even} stabilized by~PSSs \cite{IK-BJ-TG-GR-CF-MD:06}, unless sufficiently many~and carefully tuned PSSs are deployed~\cite{JHC-SGJJ-RH-WS:00,NM-LTGL:89,RAJ-PBC-NM-JCRF:10}. 
Regarding tuning of conventional PSSs, high-gain feedback is necessary in some networks~\cite{IK-BJ-TG-GR-CF-MD:06} {but may destabilize} other  networks~\cite{JHC-SGJJ-RH-WS:00,NM-LTGL:89}. {Furthermore}, even when decentralized controllers provide stability they may result in poor performance, {and their optimal tuning presents non-trivial design challenges.} 

In principle, the above problems can be solved by distributed {\em wide-area control\/} (WAC), where controllers make use of remote measurements and control signals. 
%
%
WAC is nowadays feasible thanks to recent technological advances {including} fast and reliable communication networks, high-bandwidth and time-stamped phasor measurement units (PMUs), 
and flexible AC transmission system (FACTS) devices. 
We refer to the surveys \cite{KS-SNS-SCS:06,JX-FW-CYC-KPW:06,YC-LV-GM:11} and the articles in \cite{MA:05} 
for a detailed account of technological advances and capabilities. 
%
%
{Several} efforts have been directed towards WAC of oscillations 
 based on robust and optimal control methods; see 
\cite{MA:05,QL-VV-NE:06,KT-DEB-VV-BA:05,AC:12,MZ-LM-PK-CR-GA:05,GEB-WS-JHC-TGN-NM:00,YZ-AB:08,BC-RM-BCP:04} 
and references therein.
The chosen performance metrics include frequency domain and root-locus criteria, signal amplifications from disturbance inputs to tie line flows, inter-area angles, or machine~speeds.

{%
In comparison to local control, WAC involves the communication of remote signals and may suffer from additional vulnerabilities, such as uncertain communication channels, time delays, and the lack of   globally known models and their time constants. 
{Known} time delays can be incorporated into the control design, for example, via linear matrix inequalities \cite{JHC-SGJJ-RH-WS:00,GEB-WS-JHC-TGN-NM:00} or predictor-based  design \cite{BC-RM-BCP:04}. Unknown time constants can be treated using linear parameter-varying systems \cite{QL-VV-NE:06} or online identification of parameters \cite{MZ-LM-PK-CR-GA:05}, aggregated  models \cite{AC-JHC-AS:11}, or non-parameterized models \cite{PZ-DYY-KWC-GWC:12,RE-LS:11,YP-BC-TCG:10}. Unknown delays, communication uncertainties, and uncertain time constants can be incorporated into a robust control synthesis using (un)structured uncertainties and robustness margins \cite{SS-IP:05}.}
%

Typically, the controllers are designed for {pre-specified} sensor and actuator locations and an {\em a priori\/} {fixed} sparsity pattern {that induces the} necessary communication structure. 
In an attempt to identify optimal sensor or actuator placements and to reduce the communication complexity and the interaction between control loops, 
different strategies aim at identifying few but critical control channels \cite{YZ-AB:08,AH-KI:02,NDH-LD-AFO-IK:10,LPK-RS-BCP:12,BEE-DJH:92}{; see \cite{AH-IK:08,YC-LV-GM:11} for a comparison and robustness evaluation.}
These strategies rely on modal perspectives and aim at maximizing geometric metrics such as modal controllability and observability residues. As a result, the control channels are typically chosen through combinatorial SISO criteria and not in an optimal way.

 
 Another {body} of literature relevant to our {study} is optimal control subject to structural constraints, for example, a desired sparsity pattern of the feedback matrix in~static~state feedback {design~\cite{FL-MF-MRJ:11al}}. In general, control design subject to structural constraints is hard, stabilizability is not guaranteed, and {optimal control formulations} are not convex for arbitrary structural constraints \cite{AM-NCM-MCR-SY:12}. {Furthermore}, in the absence of pre-specified structural constraints, most {optimal control formulations result in controllers without any sparsity structure according to which measurements and control channels can be selected, and they typically require a centralized implementation.}
	In order to overcome these {limitations} of decentralized optimal control, alternative strategies {have recently been proposed that simultaneously {identify} the control structure and {optimize} the closed-loop~performance; see \cite{SS-UM-FA:13,SM-UM-FA:12,UM-MP-PW:13,MF-FL-MJ:ACC11,FL-MF-MRJ:ACC12,MF-FL-MRJ:13-updated}. The proposed strategies combine classic optimal and robust control formulations with recent advances in compressed sensing {\cite{DDL:06}}.
		
Here we investigate a novel approach to WAC design for bulk power systems. We follow the sparsity-promoting optimal control approach developed in \cite{MF-FL-MRJ:13-updated} and find a linear static state feedback that {simultaneously optimizes a standard quadratic cost criterion and {induces}} a sparse control structure. 
%
%
Our choice of performance criterion is inspired by the classic work \cite{JHC:82,JHC-PK:85,DR-FD-FB:12q-updated} on slow coherency. 
We propose a novel  criterion {that encourages the closed-loop system to imitate} a homogeneous network of identical generators with no inter-area oscillations. 
	%
	%
Besides this physical insight, an additional advantage of our performance criteria is that the optimal controller makes~use of readily accessible state variables such as angles and frequencies. In the {performance} {and robustness} assessment of our sparsity-promoting {optimal} control synthesis, we emphasize robustness to time delays,  {gain} uncertainties, and variable operating conditions arising in WAC implementation.

We consider a coordinated PSS design for the IEEE 39 New England power grid model to illustrate the utility of our approach. This compelling example {shows that, with only a single WAC link,} it is possible to achieve nearly the same performance as a centralized optimal controller. In particular, our {sparsity-promoting} controller is within 1.5882\% of the performance achieved by {the optimal} centralized controller. Besides this nominal performance level, we also show that the identified control structure is not sensitive to the chosen operating point, and {that robustness margins gracefully deteriorate} as the control structure becomes increasingly sparse. 
{We also illustrate the robustness of our WAC strategy through simulations with communication noise and delays.}
Finally, we {show that our approach can also be used to identify a sparse architecture for control design via alternative means. In particular, we use nonlinear simulations to demonstrate that the closed-loop performance can be significantly improved with a proportional feedback of a single remote measurement (identified by our sparsity-promoting framework).} 


The remainder of the paper is organized as follows: Section \ref{Section: Problem Setup and Theoretical Framework} describes the sparsity-promoting optimal control approach with the slow coherency objectives. Section \ref{Section: Coordinated Supplementary PSSs Design} shows the application to the New England Power Grid. The algorithmic implementation and a link to our source code can be found in Appendix \ref{Subsection: Algorithmic Methods}. Finally, Section \ref{Section: Conclusions} concludes the paper. 

\section{Problem Setup and Theoretical Framework}
\label{Section: Problem Setup and Theoretical Framework}

\subsection{Modeling of generation, transmission, and control}


{\color{black}
We initially consider a  detailed, nonlinear, and differential-algebraic power network model of the form 
\begin{subequations}%
\begin{align}%
		\dot x(t) &= f(x(t),z(t),u(t),\eta(t)) \,,
		\label{eq: non-linear DAE model - dynamics}\\
		0 &= g(x(t),z(t),u(t),\eta(t)) \,,
		\label{eq: non-linear DAE model - constraints}
\end{align}%
\label{eq: non-linear DAE model}%
\end{subequations}%
where the dynamic and algebraic variables $x(t) \in \real^{n}$ and $z(t) \in \real^{s}$ constitute the state, $u(t) \in \real^{p}$ is the control action, and $\eta(t) \in \real^{q}$ is a white noise signal.\
Here, the dynamic equations \eqref{eq: non-linear DAE model - dynamics} account for the electromechanical dynamics of synchronous generators and their excitation control equipment. The algebraic equations \eqref{eq: non-linear DAE model - constraints} account for load flow, generator stator, and power electronic circuit equations. The control input $u(t) \in \real^{p}$ enters through either power electronics (FACTS), generator excitation (PSS) or governor control. Finally, the noise can arise from disturbances in the control loops, or fluctuating loads and generator mechanical power inputs.
Though not explicitly considered in this article, the differential-algebraic model \eqref{eq: non-linear DAE model} and our analysis are sufficiently general to allow for dynamic loads, synchronous condensers, and detailed models of FACTS devices.
}

Next we linearize the system \eqref{eq: non-linear DAE model} at a stationary operating point, solve the resulting linear algebraic equations for the variable $z(t)$, and  arrive at the linear state-space model 
\begin{equation}
	\dot x(t) = Ax(t) + B_{1} \eta(t) + B_{2} u(t)
	\label{eq: power network model}
	\,,
\end{equation}
where $A \in \real^{n \times n}$, $B_{1} \in \real^{n \times q}$, and $B_{2} \in \real^{n \times p}$. 
{The linear model \eqref{eq: power network model} is of large scale and its time constants depend on the current operating point. Thus, the precise model~\eqref{eq: power network model} of a bulk power system may not be known system-wide, or it may be known only with a limited accuracy. We discuss related robustness issues in Subsection~\ref{Subsection: Robustness and time delays}. As discussed in Section~\ref{Section: Introduction}, the model \eqref{eq: power network model} can also be estimated online 
using subspace identification methods \cite{PZ-DYY-KWC-GWC:12,RE-LS:11,YP-BC-TCG:10}, identification of aggregated models \cite{AC-JHC-AS:11}, or adaptive Kalman filtering methods \cite{MZ-LM-PK-CR-GA:05}.}

\subsection{Review of slow coherency theory} 
\label{Subsection: Review of slow coherency theory} 

To obtain an insightful perspective on inter-area oscillations, we briefly recall the classic slow coherency theory \cite{JHC:82,JHC-PK:85,DR-FD-FB:12q-updated}.
Let the state variable $x$ of the power network model \eqref{eq: non-linear DAE model}  (or its linearization \eqref{eq: power network model}) be partitioned as $x = [\theta^{T} \,,\, \dot \theta^{T} \,,\, \subscr{x}{rem}^{T}]$, where $\theta,\dot \theta \in \real^{\subscr{n}{g}}$ are the rotor angles and frequencies of $\subscr{n}{g}$ synchronous generators and $\subscr{x}{rem} \in \real^{n-2\subscr{n}{g}}$ are the remaining state variables, which typically correspond to fast electrical dynamics.
In the absence of higher-order generator dynamics, and for constant-current {and constant-impedance} loads, the power system dynamics \eqref{eq: non-linear DAE model} 
can be reduced to the electromechanical swing dynamics of the generators \cite{PK:94}:
\begin{equation}
	M_{i} \ddot \theta_{i} + D_{i} \dot \theta_{i} \!=\! \subscr{P}{red,$i$} - \sum\limits_{j=1}^{\subscr{n}{g}} |\subscr{Y}{red,$ij$}|\, E_{i} E_{j} \sin(\theta_{i} - \theta_{j} {- \varphi_{ij}})
	\,.
	\label{eq: swing eqns}
\end{equation}
Here, 
$M_{i}$ and $D_{i}$ are the inertia and damping coefficients, $E_{i}$ is the $q$-axis voltage, and $\subscr{P}{red,$i$}$ is the generator power injection in the network-reduced model. {The Kron-reduced admittance matrix $\subscr{Y}{red}$ describes the interactions among the generators, and the {phase shifts} $\varphi_{ij} \triangleq -\arctan(\Re(Y_{ij})/\Im(Y_{ij})) \in {[0,\pi/2[}$ are due to transfer conductances $\Re(Y_{ij})$ \cite{FD-FB:11d}.}
When linearized at an operating point $(\dot\theta^{*},\theta^{*})$,  the swing equations \eqref{eq: swing eqns} read as
\begin{equation}
	M \ddot \theta + D \dot \theta + L \theta = 0
	\,,
	\label{eq: linear swing eqns}
\end{equation}
where $M$ and $D$ are the diagonal matrices of\,inertia and damping coefficients, and $L$ is a Laplacian matrix with off-diagonals $L_{ij} \!=\! - |\subscr{Y}{red,$ij$}| E_{i} E_{j} \cos(\theta_{i}^{*} - \theta_{j}^{*} {-\varphi_{ij}})$ and diagonal elements $L_{ii} = - \sum_{j=1,j\neq i}^{\subscr{n}{g}} L_{ij}$. 
Notice that, {in the absence of transfer conductances}, equation \eqref{eq: linear swing eqns} describes a dissipative mechanical system with kinetic energy $(1/2) \cdot \dot \theta{^T} M \dot \theta$ and potential energy $(1/2) \cdot \theta^{T} L \theta$. Since $M$ and $D$ are diagonal, the interactions among generators in \eqref{eq: linear swing eqns} are entirely described by the weighted graph induced by the Laplacian matrix $L$.

Inter-area oscillations {may} arise from non-uniform inertia and damping coefficients (resulting in slow and fast responses), clustered groups of machines (swinging coherently) and sparse interconnections among them, as well as large inter-area power transfers.%
\footnote{A large power transfer between two nodes $i,j$ in distinct areas amounts to a large steady state difference angle $|\theta^{*}_{i} - \theta^{*}_{j}|$. Equivalently, in the linearized model  \eqref{eq: linear swing eqns}, the coupling $|L_{ij}| \!=\! |\subscr{Y}{red,$ij$}| E_{i} E_{j} \cos(\theta_{i}^{*} - \theta_{j}^{*}-{\varphi_{ij}})$ is small.
}
Assume that the set of generators {$\mc V = \left\{ 1, \ldots, n_g \right\}$} is partitioned in multiple coherent (and disjoint) groups of machines (or areas), that is, $\mc V = \mc V_{\alpha} \cup \mc V_{\beta} \cup \dots$\,. Then it can be shown \cite{JHC:82,JHC-PK:85,DR-FD-FB:12q-updated} that the long-time dynamics of each area $\alpha$ with nodal set $\mc V_{\alpha}$ are captured by the aggregate variable $\delta_{\alpha} = \bigl( \sum_{i \in \mc V_{\alpha}} M_{i} \theta_{i} \bigr) / \big( \sum_{i \in \mc V_{\alpha}} M_{i} \bigr)$ describing the center of mass of area $\alpha$. The slow inter-area dynamics are obtained as
\begin{equation}
	\tilde M \ddot \delta + \tilde D \dot \delta + \tilde L \delta = 0 \,,
	\label{eq: linear swing eqns - slow coherency}
\end{equation}
where $\delta = [\, \delta_{\alpha}~~\delta_{\beta}~\cdots ~]^{T}$  and $\tilde M$, $\tilde D$, and $\tilde L$ are the aggregated inertia, dissipation, and Laplacian matrices.

\subsection{Local and wide-area control design}

{As illustrated in the block diagram}%
\footnote{In Fig.~\ref{Fig: two-level control strategy}, the system dynamics ($A$ matrix) includes the power network dynamics and the local control loops. The WAC signal is added to the PSS input. Alternatively, the WAC signal can be added to the  automatic voltage regulator (AVR). Both strategies have been advocated, see \cite{JHC-SGJJ-RH-WS:00,GEB-WS-JHC-TGN-NM:00,YZ-AB:08}.} 
in Fig.~\ref{Fig: two-level control strategy}, we seek for linear time-invariant control laws and follow a two-level control strategy of the form  	
$u(t) = \subscr{u}{loc}(t)  + \subscr{u}{wac}(t)$.
\begin{figure}[t]
	\centering{
	\includegraphics[width=0.95\columnwidth]{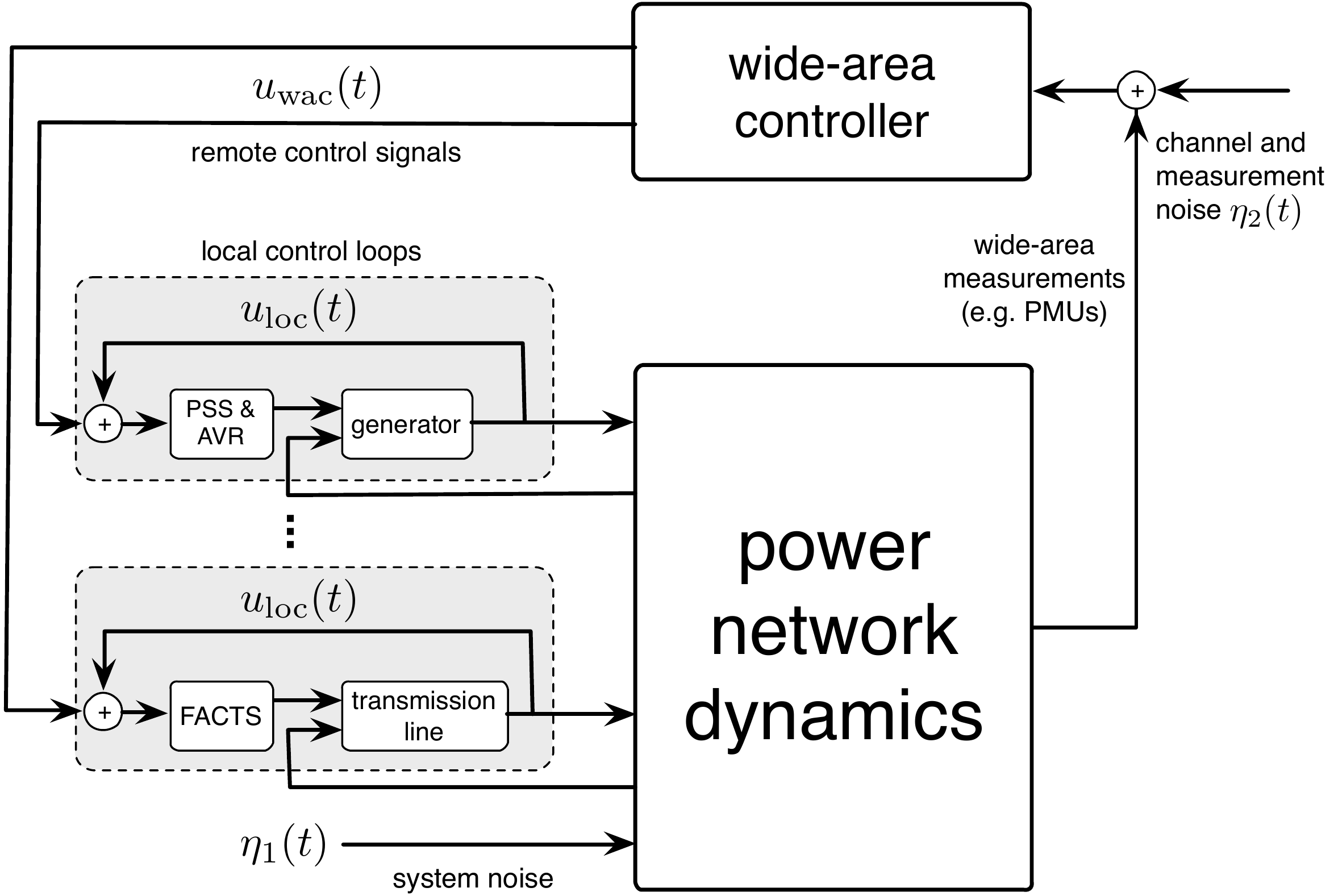}
	\caption{Two-level control design combining local and wide-area control.}
	\label{Fig: two-level control strategy}
	}
\end{figure}

In the first step, the local control $\subscr{u}{loc}(t)$ is designed based on locally available measurements and with the objective of stabilizing each isolated component. For example, $\subscr{u}{loc}(t)$ can be obtained by a conventional PSS design with the objective {of suppressing} local oscillations~\cite{PK:94}. Next, the wide-area control $\subscr{u}{wac}(t)$ is designed {to enhance} the global system behavior and {suppress} inter-area oscillations. 
For this second design step, the local control $\subscr{u}{loc}(t)$ is assumed to be absorbed into the plant \eqref{eq: power network model}.
Additionally, since $\subscr{u}{wac}(t)$ relies on the communication of remote signals, the two-level control strategy guarantees a nominal performance level in case of communication~failures.
%

\subsection{Sparsity-promoting linear quadratic control}

As discussed in Section \ref{Section: Introduction}, an inherent problem in WAC is the proper choice of {\em control architecture}, that specifies which quantities need to be measured and which controller needs to access which measurement.
Here, we confine our attention to static state feedback   
$
	\subscr{u}{wac}(t) 
	~=\; - K \, x(t),
$
where the control structure is determined by the sparsity pattern of the feedback gain {} $K \in \real^{p \times n}$. 
%
%
We use the sparsity-promoting optimal control framework~\cite{MF-FL-MRJ:13-updated} to minimize the $\ell_1$-regularized steady-state variance of a stochastically-driven closed-loop system:
	\begin{equation}
	\!\!\!
	\begin{array}{rl}
\minimize & \displaystyle{\lim_{t \, \to \, \infty}} \mc E \left\{ x(t)^{T} Qx(t) + u(t)^{T}Ru(t) \right\} 
	\\[0.15cm]
&\; + ~\Bigl. \gamma \; \displaystyle{\sum\nolimits_{i,j}} w_{ij}|K_{ij}|  \\
\textup{subject to} & \\[0.15cm]
\textup{dynamics:} & \dot x(t) = Ax(t) + B_{2}u(t) + B_{1} \eta(t),\\[0.15cm]
\textup{linear control:} & u(t) = - Kx(t), \\[0.15cm]
\textup{stability:} & \bigl( A-B_{2}K \bigr) \textup{ Hurwitz}.
	\label{eq: optimal control problem}
	\end{array}
	\end{equation}
Here, $\gamma \geq 0$ is a nonnegative parameter, $\mc E \{ \cdot \}$ is the expectation operator, and $Q \in \real^{n \times n}$, $R \in \real^{n \times n}$ are positive semidefinite and positive definite matrices that denote the state and control weights, respectively. We assume that the triple $(A,B,Q^{1/2})$ is stabilizable and detectable.
The term $\sum_{i,j} w_{ij}|K_{ij}|$ is a {\em weighted $\ell_{1}$-norm\/} of $K$, where $w_{ij} >0$ are positive weights. The weighted $\ell_{1}$-norm serves as a proxy for the (non-convex) cardinality function $\textup{card}(K)$ denoting the number of non-zero entries {in} $K$. 
{An effective method for enhancing sparsity is} to solve a sequence of weighted $\ell_{1}$-optimization problems, where the weights are determined by the solution of the weighted $\ell_{1}$-problem in the previous iteration; see \cite{EJC-MBW-SPB:08} for further details.

 {An equivalent formulation of the optimal control problem \eqref{eq: optimal control problem} is given via the closed-loop observability Gramian $P$ as
	\begin{align}
		\displaystyle{\minimize}\quad
		&
		J_{\gamma}(K)  \; \triangleq  \; \textup{trace} \, \bigl(B_1^{T}  P  B_1\bigr) + \gamma  \;\displaystyle{\sum\nolimits_{i,j}} \; w_{ij} \, |K_{ij}| 
		\nonumber\\
	    \textup{subject to}\quad 
	&	    
	    \bigl(A-B_{2}K)^{T} P + P(A-B_{2}K)
	  \label{eq: optimal control problem - H2}\\&\hfill\qquad\qquad\qquad\qquad\qquad
	   = - (Q+K^{T}RK).
	\nonumber
	\end{align}
	The latter formulation \eqref{eq: optimal control problem - H2} is amenable to an iterative solution strategy using the alternating direction method of multipliers (ADMM); see \cite{FL-MF-MRJ:ACC12,MF-FL-MRJ:13-updated}.} The cost function in \eqref{eq: optimal control problem - H2} is  a linear combination of the $\mc H_{2}$-norm of the closed-loop system 
\begin{align*}
	\dot x(t) &= \left(A-B_{2}K\right) x(t) + B_{1} \eta(t), 
	\\
	{
	y(t)
	} 
	&= 
	{
	\left[
	\begin{array}{c}
	Q^{1/2} 
	\\
	- R^{1/2} K 
	\end{array}
	\right] 
	x(t),
	}
\end{align*}
and the sparsity-promoting term {$\gamma \sum_{i,j} w_{ij}|K_{ij}|$}. 
In what follows, for a fixed value of $\gamma \geq 0$, we denote the minimizer to \eqref{eq: optimal control problem - H2} by $K^{*}_{\gamma}$ and the {optimal} cost by $J^{*}_{\gamma} = J(K^{*}_{\gamma})$. For $\gamma = 0$ the problem \eqref{eq: optimal control problem - H2} reduces to the {standard state-feedback $\mc H_2$-problem \cite{SS-IP:05} with the optimal gain $K^{*}_{0}$ and the optimal cost $J_{0}^{*}$}. On the other hand, for $\gamma > 0$ the {weighted $\ell_1$-norm promotes sparsity in the feedback gain $K^{*}_{\gamma}$, thereby identifying essential pairs of control inputs and measured outputs.}

{To solve the sparsity-promoting optimal control problem\,\eqref{eq: optimal control problem - H2}, we rely on the approach proposed in~\cite{MF-FL-MRJ:13-updated} 
using ADMM.
See Appendix~\ref{Subsection: Algorithmic Methods} for details on the algorithmic implementation.}


\subsection{Choice of optimization objectives}

The design parameters $Q$, $R$, $B_{1}$, $\gamma$ need to be chosen with the objective of damping inter-area oscillations. The resulting feedback $\subscr{u}{wac}(t)$, the control variables $K_{\gamma}^{*}x(t)$, the communication structure ({the sparsity pattern of the off-diagonals of $K_{\gamma}^{*}$), and the control effort depend solely on {the system model and our choices of the design parameters} $Q$, $R$, $B_{1}$, and\,$\gamma$.

{\bf\em State cost:} Slow coherency theory shows that an {\em ideal} power system without inter-area oscillations is characterized by uniform inertia coefficients and homogenous power transfers among generators. Equivalently, the Laplacian and inertia matrices in the linearized swing equations \eqref{eq: linear swing eqns} take the form
\begin{equation*}
	L = \subscr{L}{unif} = \ell \cdot \bigl( I_{\subscr{n}{g}} - (1/\subscr{n}{g}) \mathds{1}_{\subscr{n}{g}} \mathds{1}_{\subscr{n}{g}}^{T} \bigr)
	\;,~~ 
	M = \subscr{M}{unif} = m \cdot I_{\subscr{n}{g}}
	\,,
\end{equation*}
where $\ell,m > 0$ are constants, $I_{\subscr{n}{g}}$ is the $\subscr{n}{g}$-dimensional identity matrix, and $\mathds{1}_{\subscr{n}{g}}$ is the $\subscr{n}{g}$-dimensional vector of ones.

Inspired by the above considerations, we choose the following performance specifications for the state cost: 
\begin{equation}
	x^{T} Q x = \frac{1}{2} \theta^{T} \subscr{L}{unif}\, \theta + \frac{1}{2} \dot \theta^{T} \subscr{M}{unif}\, \dot \theta + \varepsilon \cdot \| \theta \|_{2}^{2} 
	\label{eq: state cost -- average performance}
	\,,
\end{equation}
where $\varepsilon \geq 0$ is a small regularization parameter. 
For $\varepsilon = 0$, the state cost $x^{T}Qx$ {quantifies} the kinetic and potential energy of a homogenous network composed of identical generators, and it penalizes frequency violations and angular differences. 
The cost $x^{T}Qx$ does not penalize deviations in the generator voltages, the states of the excitation system, or the local~control $\subscr{u}{loc}(t)$ included in the $A$ matrix.
The regularization term $\varepsilon \| \theta \|_{2}^{2}$ assures  numerical stability and detectability%
\footnote{If the power system model \eqref{eq: non-linear DAE model} entails a slack bus, then $\varepsilon$ can be set to zero. Otherwise, there is a marginally stable and unobservable mode corresponding to the {rotational symmetry of the power flow: a uniform shift of all angles does not change the vector field \eqref{eq: non-linear DAE model}.} 
The optimal control problem \eqref{eq: optimal control problem - H2}~is~still feasible, and the marginally stable mode is not affected by controllers using only angular differences (i.e., power flows) \cite{BB-FDG:12-updated,MF-FL-MRJ:13-updated}. To guarantee a simple numerical treatment of the resulting Ricatti equations using standard software, an inclusion of a small positive regularization parameter $\varepsilon>0$ is necessary.} of $(A,Q^{1/2})$.

As we will see in Section \ref{Section: Coordinated Supplementary PSSs Design}, the state cost \eqref{eq: state cost -- average performance} results~in an improved {\em average} closed-loop performance with all inter-area modes either damped or distorted.
If the objective~is~to reject a {\em specific} inter-area mode, 
for example, a dominant inter-area mode featuring the two groups $\mc V_{\alpha}$ and $\mc V_{\beta}$, then the discussion preceding the inter-area dynamics \eqref{eq: linear swing eqns - slow coherency} suggests the~cost 
\begin{align}
	x^{T} Q x =&\;\;
	\ell \cdot \bigl\| \delta_{\alpha} - \delta_{\beta} \bigr\|_{2}^{2} +  m \cdot \bigl\| \dot\delta_{\alpha} - \dot\delta_{\beta} \bigr\|_{2}^{2} 
	+ \varepsilon \cdot \| \theta \|_{2}^{2} 
	\label{eq: state cost -- particular mode}
	\,,
\end{align}
where $\ell,m > 0$ are gains, $\delta_{\alpha},\delta_{\beta}$ are the aggregate variables, and $\varepsilon\geq0$ is a small regularization parameter.
Alternative cost functions penalize selected generator frequency deviations or branch power flows $|L_{ij}| (\theta_{i} - \theta_{j})$ to assure coherency and guarantee (soft) thermal limit constraints. Finally, linear combinations of all cost functions can also be chosen.
 
{In the costs \eqref{eq: state cost -- average performance} and \eqref{eq: state cost -- particular mode}, the parameter $m$ mainly affects the frequency damping, whereas the parameter $\ell$ affects the damping of the machine difference angles, that is, inter-area power flows.} In summary, the state costs \eqref{eq: state cost -- average performance}-\eqref{eq: state cost -- particular mode} reflect the insights of slow coherency theory, 
and, as we will see later, they also promote the use of readily available control variables.


{\bf\em Control cost:} For simplicity and in order to minimize interactions among generators the control effort is penalized as $u^{T}Ru$, where $R$ is a positive definite  diagonal matrix. {Larger diagonal values of $R$ result in a smaller control effort.}

{\bf\em System noise:}
{In standard (centralized) linear quadratic control, the optimal feedback gain $K_{0}^{*}$ is {independent} of the choice of $B_{1}$ \cite{SS-IP:05}. However, for $\gamma>0$ the optimal feedback gain $K_{\gamma}^{*}$ depends on the choice of $B_{1}$; see the necessary optimality conditions in \cite{MF-FL-MRJ:13-updated}.} In order to mitigate the impacts of noisy or lossy communication among spatially distributed controllers, one may choose $B_{1} = B_{2}$. Otherwise, $B_{1}$ can be chosen to include load and generation uncertainties in~\eqref{eq: non-linear DAE model}.

{\bf\em Promoting sparsity:}
For $\gamma = 0$, the optimal control problem \eqref{eq: optimal control problem - H2} reduces to a {standard state-feedback $\mc H_{2}$}-problem whose solution can be obtained from the {positive definite} solution to the algebraic Riccati equation. Starting from this initial value, we iteratively solve the {optimal control problem} \eqref{eq: optimal control problem - H2} for increasingly larger values of $\gamma$. {We found that a logarithmically spaced sequence of {$\gamma$-}values performs well in practice.} In the end, the resulting sequence of optimal controllers is analyzed, and a value of $\gamma$ is chosen {to strike} a balance between the closed-loop performance and sparsity of the controller.

{In summary, if the state costs \eqref{eq: state cost -- average performance} or \eqref{eq: state cost -- particular mode} are chosen, then the parameters affecting the  optimal feedback gains $K_{\gamma}^{*}$ are $\ell$, $m$, $\varepsilon$, $B_{1}$, the diagonal values of $R$, and the sequence of $\gamma$-values.}

\subsection{Robustness, time delays, and {gain} uncertainties}
\label{Subsection: Robustness and time delays}

\begin{figure}[t]
	\centering{
	\includegraphics[width=0.89\columnwidth]{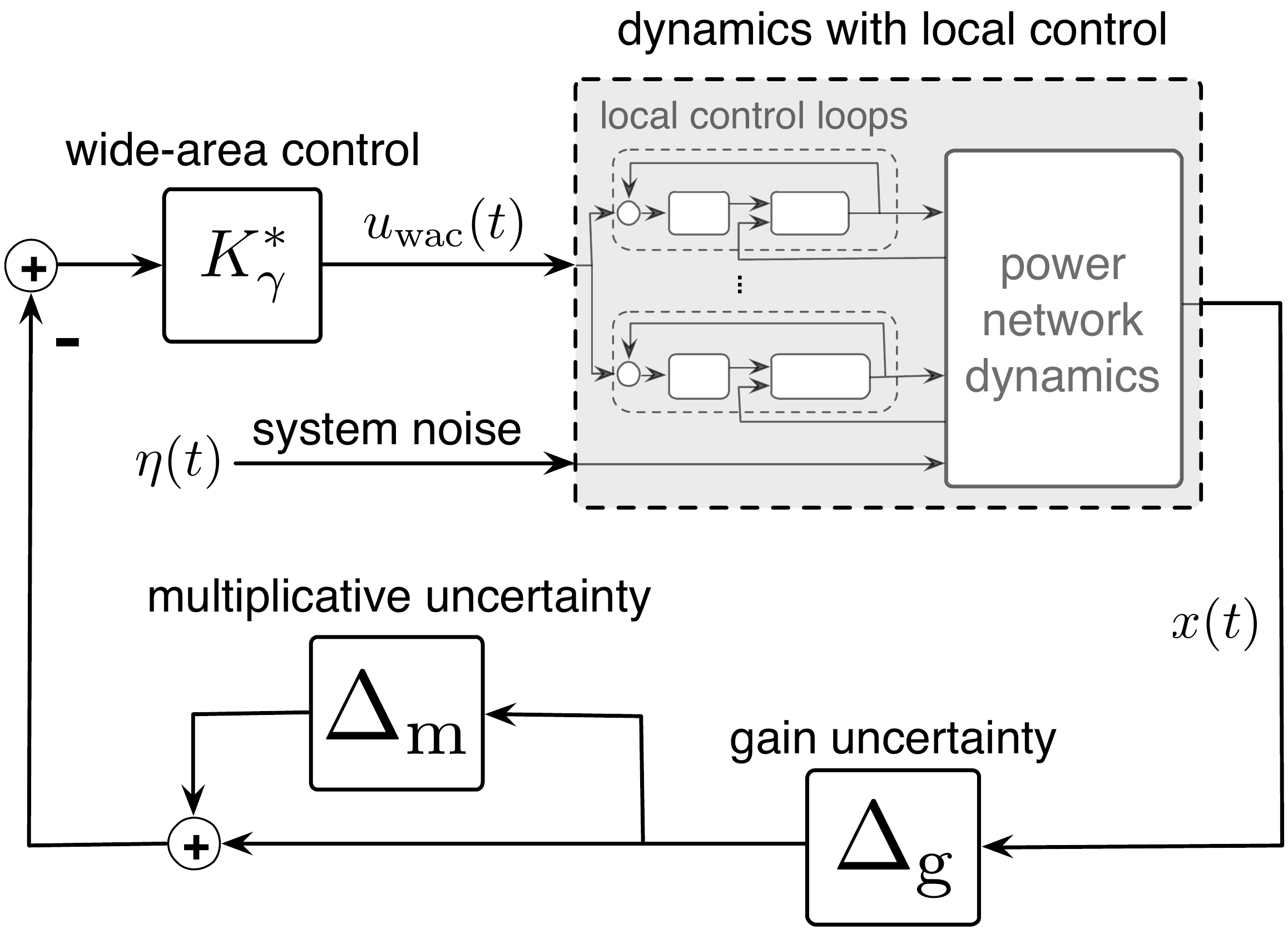}
	\caption{Gain uncertainties $\subscr{\Delta}{g}$ and multiplicative uncertainties $\subscr{\Delta}{m}$ account for delays and uncertain (or unmodeled) dynamics in the WAC channel.}
	\label{Fig: uncertainty loop}
	}
\end{figure}


{As discussed in Section \ref{Section: Introduction}, WAC may have to deal with uncertain communication channels, time delays, and the lack of  precisely known system-wide time constants. In particular, time delays {may} arise from communication delays, latencies and multiple data rates in the SCADA (supervisory control and data acquisition) network, asynchronous measurements, and local processing and operating times. As a result, local control signals and measurements may have different rates and time stamps than WAC signals obtained from remote sites.}

Given the above vulnerabilities, the wide-area control loop needs to be robust to gain uncertainties $\subscr{\Delta}{g}$ and multiplicative uncertainties $\subscr{\Delta}{m}$, as depicted in Fig.~\ref{Fig: uncertainty loop}. 
For instance, if the channel from measurement $i$ to control input $j$ features a time-delay, then the $(i,j)${th} entry of $\subscr{\Delta}{g}$ is a pure phase delay. We refer to \cite[Section 7.4]{SS-IP:05} for the modeling of delays with multiplicative uncertainties $\subscr{\Delta}{m}$.
In a general robust control theory framework, the uncertainty blocks $\subscr{\Delta}{g}$ and $\subscr{\Delta}{m}$ are unknown, stable, and proper dynamical systems that typically satisfy certain norm bounds on their input-output behavior. 

{Standard} frequency domain arguments show that robustness to delays is directly related to phase margins \cite{SS-IP:05}.
For a non-interacting  (i.e., diagonal) uncertainty $\subscr{\Delta}{g}$, the {optimal centralized feedback gain $K_{0}^{*}$} results in guaranteed phase margins of $\pm 60^{\circ}$, lower gain (reduction) margins of 0.5, and upper gain (amplification) margins of $\infty$ in every channel~\cite{MS-MA:77}. 
	%
	%
For an unstructured multiplicative uncertainty $\subscr{\Delta}{m}$, the optimal centralized feedback gain $K_{0}^{*}$ {guarantees closed-loop stability} provided that its $\mathcal L_{2}$-induced input-output norm is bounded by 0.5 corresponding to a multivariable phase margin of $30^{\circ}$~\cite{JD-GS:81}. 
%

{
Yet another source of uncertainty results from the fact that the time constants of the linearized model \eqref{eq: power network model} depend on the current operating condition of the overall power system. The gain margins discussed in preceding paragraphs {provide} safeguards against gain uncertainties in the control loop arising, for example, from variable or uncertain time constants in the generator excitation system or the measurement channels.}

	%
	%
For {the $\gamma$-parameterized family of sparsity-promoting} controllers, we explicitly verify robustness to phase and gain uncertainties in Section~\ref{Section: Coordinated Supplementary PSSs Design}. In particular, we observe that the phase and gain margins decay gracefully for $K_{\gamma}^{*}$ as $\gamma$ increases.
Alternatively, robustness to {\em known\/} time delays can be included in the control design by explicitly {accounting for} delays {via} Pad\'e approximations (absorbed in the~plant)~\cite{SS-IP:05}. 


\section{Coordinated Supplementary PSSs Design}
\label{Section: Coordinated Supplementary PSSs Design}

In this section, we {use the IEEE 39 New England power grid model} to {validate performance of} the proposed WAC strategy. {This model consists} of 39 buses and 10 two-axis generator models, where generator 10 is an equivalent aggregated model.

\subsection{Local control design and inter-area dynamics}

The {\em Power System Toolbox} \cite{JHC-KWC:92} was used to obtain the nonlinear differential-algebraic model \eqref{eq: non-linear DAE model} and the linear state space system~\eqref{eq: power network model}. 
%
The open-loop system is unstable, and the generators are equipped with PSS excitation controllers designed with washout filters and lead/lag elements. %
For generator $i$, the local PSS reads in the Laplace domain~as
\begin{equation}
	\subscr{u}{loc,$i$}(s) = k_{i} \cdot \frac{\subscr{T}{w,$i$}s}{1+\subscr{T}{w,$i$}s} \cdot \frac{1+\subscr{T}{n1,$i$}s}{1+\subscr{T}{d1,$i$}s} \cdot \frac{1+\subscr{T}{n2,$i$}s}{1+\subscr{T}{d2,$i$}s} \, \cdot \dot \theta_{i}(s)
	\label{eq: PSS design}
	\,,
\end{equation}
{where, with an abuse of notation, $\dot\theta_{i}(s)$ and $\subscr{u}{loc,$i$}(s)$ denote the Laplace transforms of $\dot \theta_{i}(t)$ and $\subscr{u}{loc,$i$}(t)$, respectively.}
The controller gains are chosen according to the tuning strategy \cite{RAJ-PBC-NM-JCRF:10} 
as $\subscr{T}{w,$i$} = 3$, $\subscr{T}{n1,$i$} = \subscr{T}{n2,$i$} = 0.1$, $\subscr{T}{d1,$i$} = \subscr{T}{d2,$i$} = 0.01$ for $i \in \until{9}$, $k_{i} = 12$ for $i \in \{1,2,3,5,6,9\}$, $k_{4} = 10$, $k_{7} = 11.03$, and $k_{8} = 9.51$.
%
As in \cite{RAJ-PBC-NM-JCRF:10}, we adopt the first-order exciter model, and choose uniform time constants $\subscr{T}{a} = 0.015\,\textup{s}$ and regulator gains $\subscr{K}{a} = 200$ for all generators. The static DC gains $\subscr{K}{a} \cdot k_{i}$ in the exciter and PSS control channels are slightly different compared to the highly tuned data in \cite{RAJ-PBC-NM-JCRF:10}. We made this choice for the sake of illustration: the PSSs \eqref{eq: PSS design} provide good damping for the local modes and stabilize the otherwise unstable open-loop system \eqref{eq: power network model}. On the other hand, the chosen PSSs do not severely distort the shapes of the inter-area modes unless they are very carefully tuned to the linearized model (depending on the particular operating condition). 

An analysis of the closed-loop modes and participation factors reveals the presence of five dominant inter-area modes. These five modes are reported in Table \ref{Table: inter-area modes of New England Grid}, and the groups of coherent machines and  the frequency components of the associated eigenvectors are illustrated in Fig.~\ref{Fig: New England grid with inter-area modes}.

\begin{table}[h]
\centering
\caption{Inter-area modes of New England power grid with PSSs}
\label{Table: inter-area modes of New England Grid}
\begin{tabular}{| c | l | l | l | l |}
\hline
\!\!mode\!\! \!\!\!&\!\!\! eigenvalue \!\!\!&\!\!\! damping \!\!\!&\!\!\! frequency \!\!\!&\!\!\! coherent \\
no. \!\!\!&\!\!\! pair  \!\!\!&\!\!\! ratio \!\!\!&\!\!\! $[\textup{Hz}]$ \!\!\!&\!\!\! groups \\
\hline\hline
\!\!1 \!\!\!&\!\!\! $-0.6347 \pm \textup{i}\, 3.7672$ \!\!\!&\!\!\! 0.16614 \!\!\!&\!\!\! 0.59956 \!\!\!&\!\!\! 10 vs. all others \!\!\!\!\\
\!\!2 \!\!\!&\!\!\! $-0.7738 \pm \textup{i}\, 6.7684$ \!\!\!&\!\!\! 0.11358 \!\!\!&\!\!\! 1.0772 \!\!\!&\!\!\! 1,8 vs. 2-7,9,10 \!\!\!\!\\
\!\!3 \!\!\!&\!\!\! $-1.1310 \pm \textup{i}\, 5.7304$ \!\!\!&\!\!\! 0.19364 \!\!\!&\!\!\! 0.91202 \!\!\!&\!\!\! 1,2,3,8,9 vs. 4-7 \!\!\!\!\\
\!\!4 \!\!\!&\!\!\! $-1.1467 \pm \textup{i}\, 5.9095$ \!\!\!&\!\!\! 0.19049 \!\!\!&\!\!\! 0.94052 \!\!\!&\!\!\! 4,5,6,7,9 vs. 2,3 \!\!\!\!\\
\!\!5 \!\!\!&\!\!\! $-1.5219 \pm \textup{i}\, 5.8923$ \!\!\!&\!\!\! 0.25009 \!\!\!&\!\!\! 0.93778 \!\!\!&\!\!\! 4,5 vs. 6,7 \!\!\!\!\\
\hline
\end{tabular}
\end{table}

\begin{figure}[htbp]
	\centering{
\includegraphics[width=0.9\columnwidth]{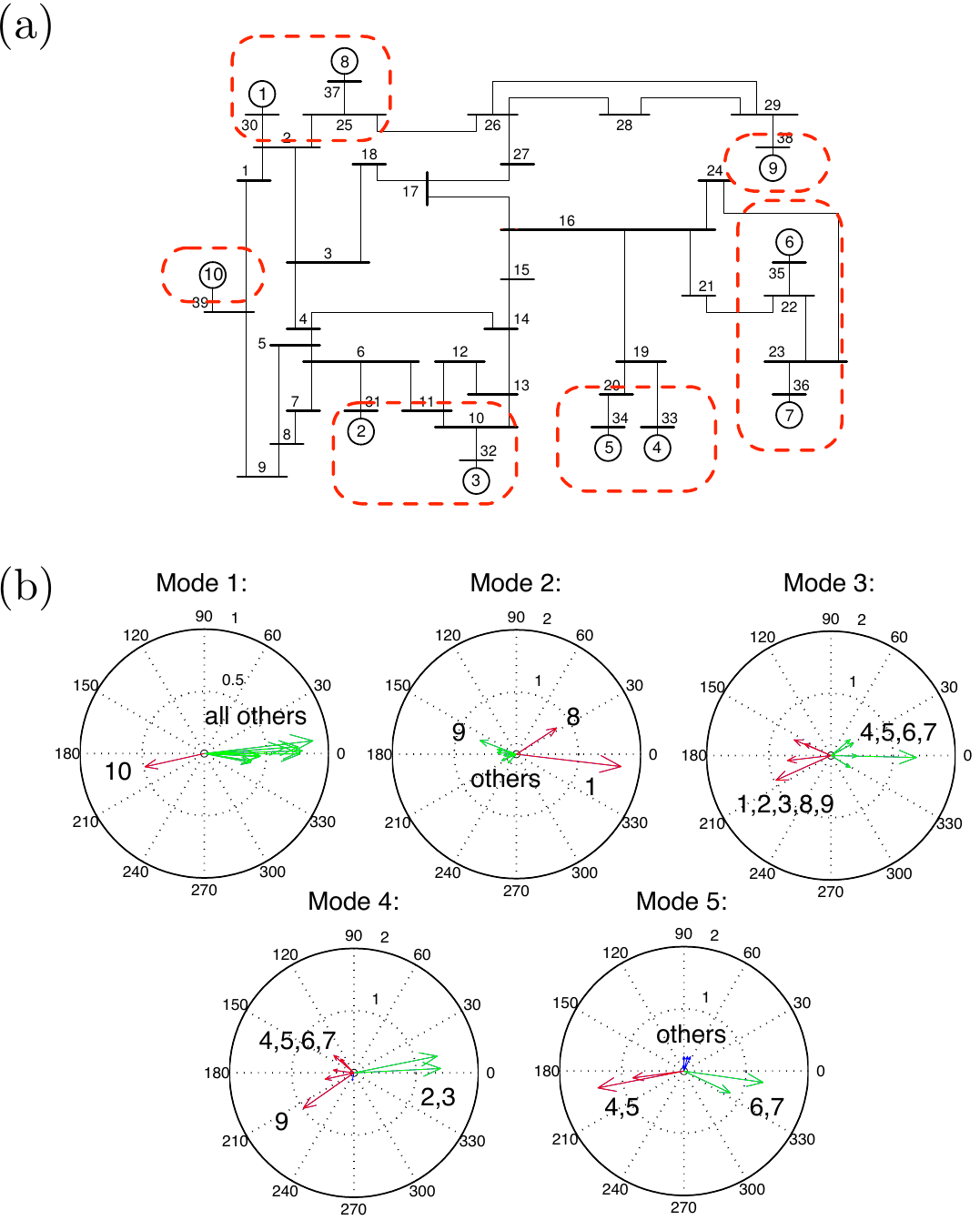}
	\caption{Subfigure (a) displays the IEEE 39 New England power grid and its coherent groups. 
	The polar plots in Subfigure (b) show the generator frequency components of the five {dominant} inter-area modes.}
	\label{Fig: New England grid with inter-area modes}
	}
\end{figure}

\subsection{WAC design and nominal performance}

To provide additional damping for the inter-area oscillations, we design a supplementary WAC signal $\subscr{u}{wac}(t)$ which additively enters the AVRs of the controlled generators. {The sparsity-promoting optimal control strategy \eqref{eq: optimal control problem} is used to design $\subscr{u}{wac}(t)$,} where the state cost \eqref{eq: state cost -- average performance} is selected with~$\varepsilon  = 0.1$ and $(\ell,m) = (2,2)$.
To share the control burden equally we set the control weight to be identity $R = I$. This choice results in a WAC signal $\subscr{u}{wac}(t)$ {of the same magnitude as the local control signal} $\subscr{u}{loc}(t)$, that is, $\max_{t \in \real_{\geq 0}}\| \subscr{u}{wac}(t) \|_{{\infty}} \approx \max_{t \in \real_{\geq 0}}\| \subscr{u}{loc}(t) \|_{{\infty}}$, and it avoids input saturation.
{Furthermore, to reject communication noise in the WAC implementation, we choose $B_{1}=B_{2}$. Finally, we solve the optimal control problem \eqref{eq: optimal control problem - H2} for 40 logarithmically spaced values of $\gamma$ in the interval $[10^{-4},10^{0}]$. Our results are reported in Fig.~\ref{Fig: NE - performance}-\ref{Fig: NE - Spy(BK)}.}

\begin{figure}[tb]
	\centering{
	\includegraphics[width=0.99\columnwidth]{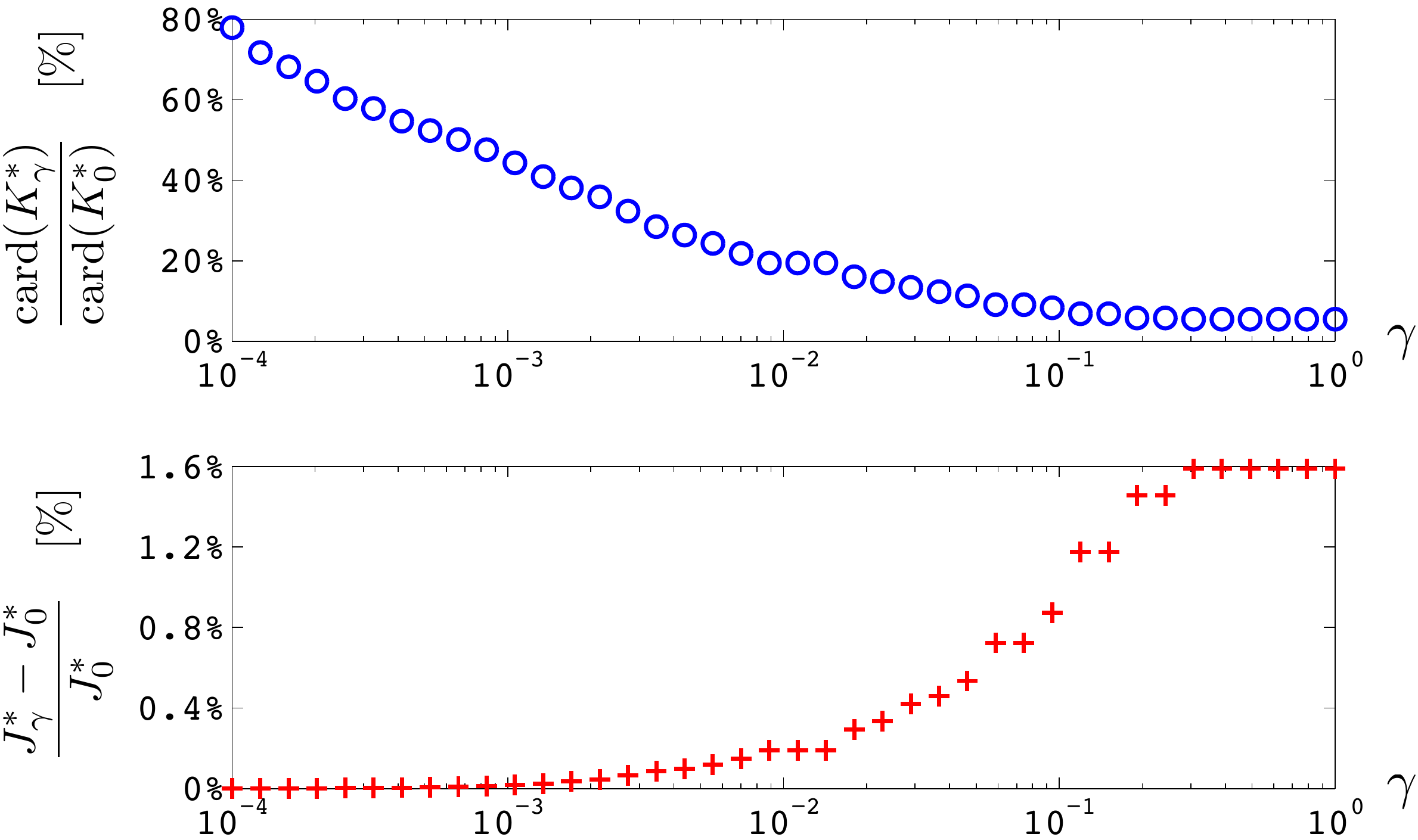}
	}
	\caption{Number of nonzero entries in {the sparse feedback gain} $K^{*}_{\gamma}$ {and the performance degradation of $K^*_\gamma$
	relative to the optimal centralized gain $K^*_0$.}}
	\label{Fig: NE - performance}
\end{figure}
\begin{figure}[tb]
	\centering{
	\includegraphics[width=0.95\columnwidth]{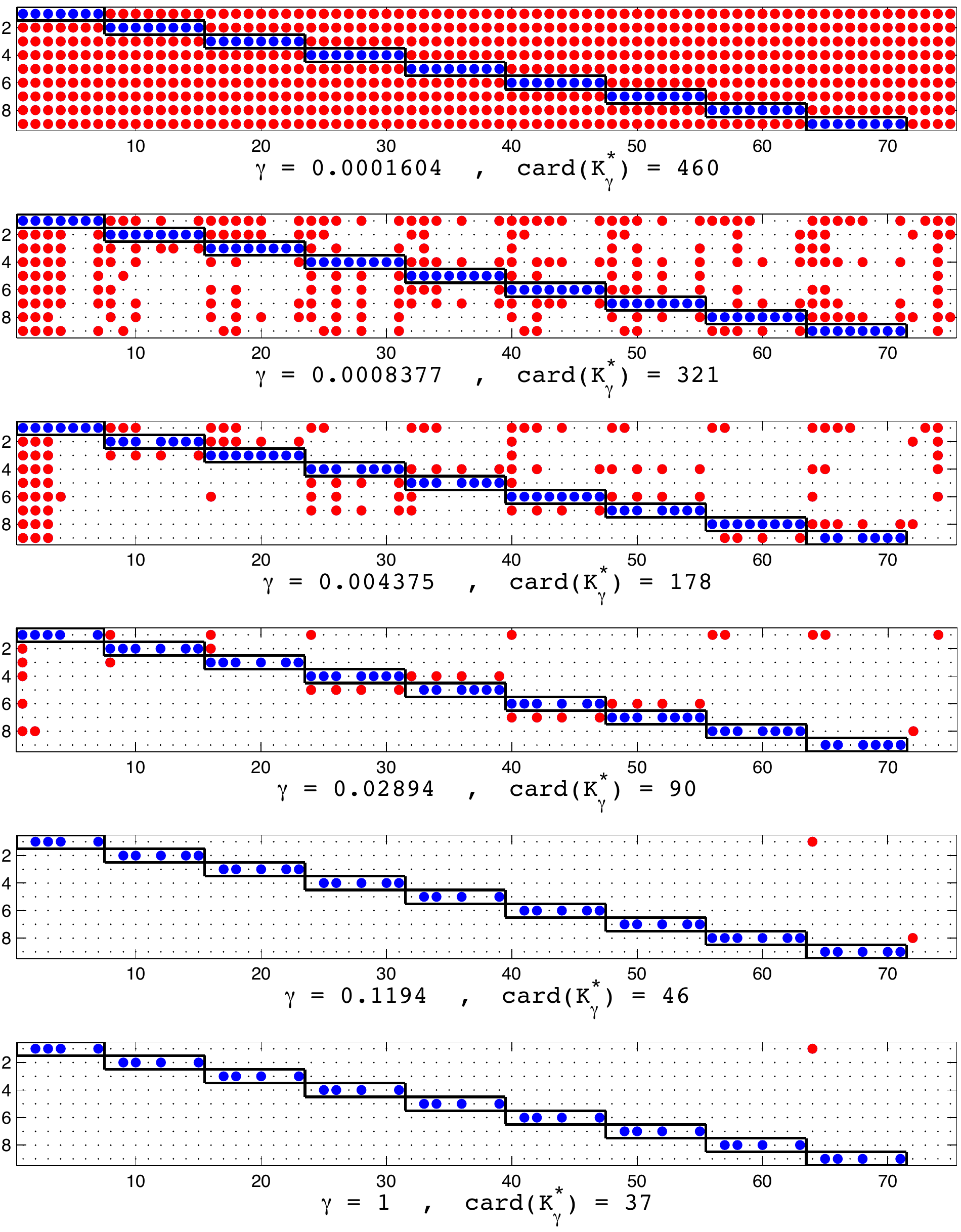}
	}
	\caption{The sparsity patterns of $K_{\gamma}^{*}$ illustrate the control architecture. A nonzero $(i,j)$ element of $K_{\gamma}^{*}$ (denoted by a colored dot) implies that controller $i$ needs to access state $j$. The diagonal blocks (framed) with blue dots correspond to local feedback, and the off-diagonal blocks with red dots correspond to remote feedback signals. As $\gamma$ increases, the information exchange becomes sparser, and angles and frequencies (the first two states of each vertical block) become the sole signals to be communicated.}
	\label{Fig: NE - Spy(BK)}
\end{figure}

%
For $\gamma = 0$, the optimal feedback gain $K^{*}_{0}$ is fully populated, {thereby requiring centralized implementation.}
As $\gamma$ increases, the off-diagonal elements of the feedback matrix $K_{\gamma}^{*}$ become significantly sparser whereas the relative cost $\bigl( J^{*}_{\gamma} - J^{*}_{0} \bigr) / J^{*}_{0} $ increases only slightly; see Fig.~\ref{Fig: NE - performance} and Fig.~\ref{Fig: NE - Spy(BK)}. 
Additionally, as $\gamma$ increases, the state cost \eqref{eq: state cost -- average performance} {promotes the use of angles and speeds in the off-diagonal elements of $K_{\gamma}^{*}$, and most nonzero elements of $K_{\gamma}^{*}$ correspond to local feedback.}
{For $\gamma = 1$, the optimal} controller $K_{1}^{*}$ is within 1.5882\% of the optimal centralized performance even though only a single signal needs to be communicated: the controller at generator 1 needs to access $\theta_{9}(t)$.
As expected, as $\gamma$ increases most of the control burden is on generator 1, which has the largest inertia of all controlled generators. Likewise, the angle of loosely connected%
\footnote{In terms of graph theory, the sum of effective resistances \cite{FD-FB:11d} between generator 9 and the others is very large compared to remaining~network.}
generator 9 needs to be measured and communicated to generator 1. The identified WAC channel $9 \to 1$ appears to be necessary to suppress the inter-area mode 2 in Table \ref{Table: inter-area modes of New England Grid}, which is mainly dominated by generators 1\,and\,9.
%

\begin{figure}[htbp]
	\centering{
	\includegraphics[width=0.99\columnwidth]{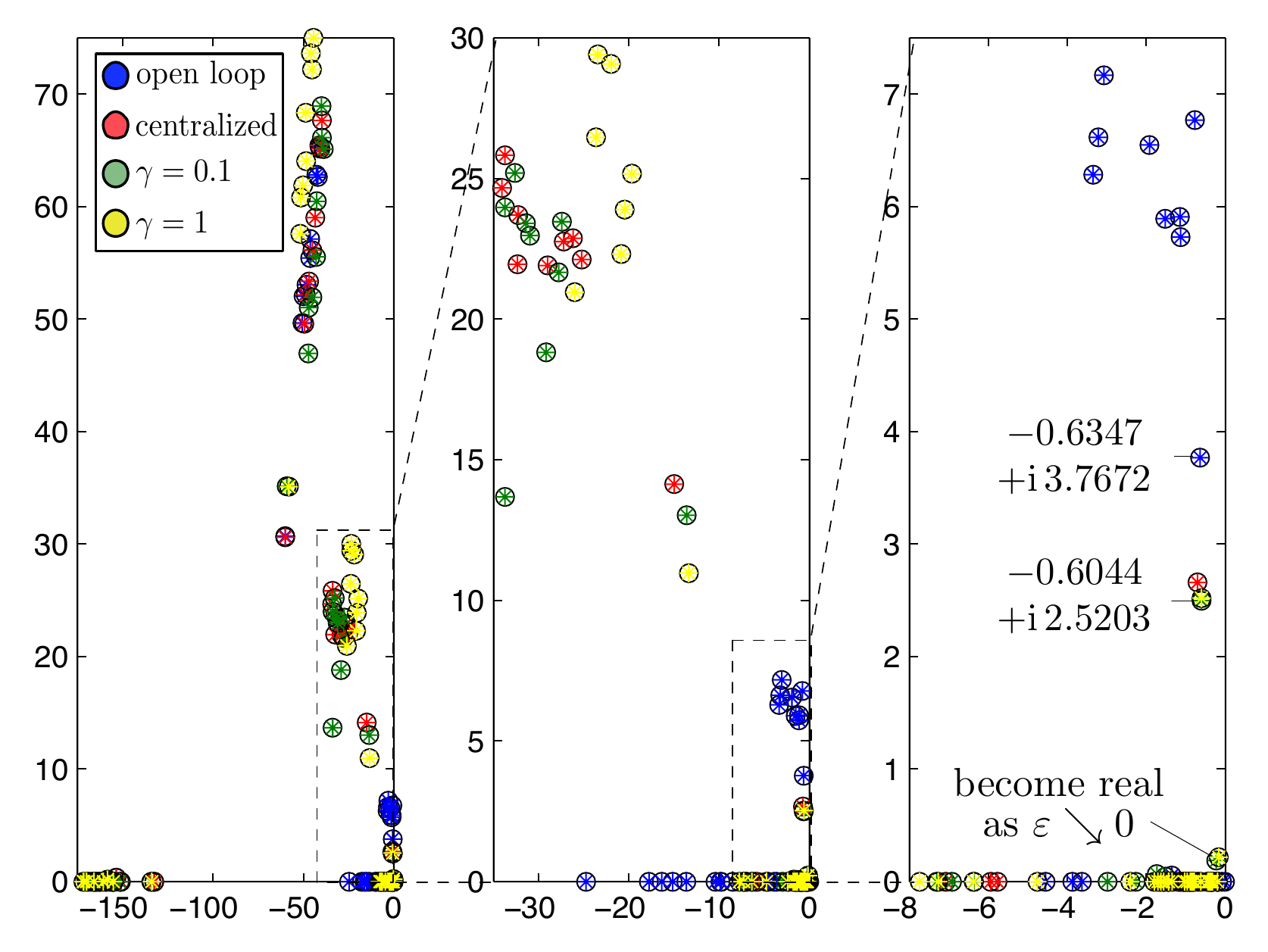}
	\caption{The spectrum of the WAC open loop matrix $A$ displayed in blue and the spectrum of the WAC closed-loop matrix $A-B_{2}K_{\gamma}$ for $\gamma \in \{0,0.1,1\}$ displayed in red, green, and yellow, respectively. The central and right panels are magnifications of the left and central panel, respectively.}
	\label{Fig: spectrum}
	}
\end{figure}
{Fig.~\ref{Fig: spectrum} displays the spectrum of the open-loop matrix $A$ and the closed-loop matrix $A-B_{2}K_{\gamma}$ for different values of $\gamma$. We make the following observations: First, as $\gamma$ increases, the locations of the poles of $A-B_{2}K_{\gamma}$  do not change significantly, and the poorly damped poles essentially overlap. This confirms our observation from Fig.~\ref{Fig: NE - performance}:  ${K}_{\gamma}$ (with $\gamma >0$) achieves a similar performance as the centralized feedback gain $K_{0}$. Second, the weakly damped pair of poles (with with nearly zero imaginary part and very high damping ratio $0.7452$) is an artifact of the $\varepsilon$-regularization in \eqref{eq: state cost -- average performance}. For $\varepsilon \searrow 0$, this pole pair becomes real and one of the poles becomes zero corresponding to the rotational symmetry of the power flow. Third, in comparison with the open loop, the wide-area controller does not significantly change the real part of the eigenvalues of the weakly damped oscillatory modes. Rather, the {closed-loop performance is improved by increasing the damping ratio and distorting the} associated eigenvectors. Apart from one pair of eigenvalues located at $-0.6044 \pm\textup{i}\,  2.5203 $ with high damping ratio $0.239$ (possibly corresponding to the inter-area mode 1 in Table~\ref{Table: inter-area modes of New England Grid}),} all {other} complex-conjugate eigenvalue pairs are left of the asymptote $\textup{Real}(s) = -12.74$. From the frequency components of its  eigenvector (in Fig.~\ref{Fig: time-series NE 1}(a)) and the {frequency time series} (in {Fig.~\ref{Fig: time-series NE 2}(a) and \ref{Fig: time-series NE 2}(b))}, it can be seen that this mode does not anymore correspond to generators oscillating against each other. As a consequence, poorly damped power flow oscillations between the areas $\mc V_{\alpha} \!=\! \until 9$ and $\mc V_{\beta} \!=\! \{10\}$ {have been suppressed}; see {Fig.~\ref{Fig: time-series NE 2}(d) and~\ref{Fig: time-series NE 2}(e).} {We finally note that, besides the critical mode, the damping ratio is also significantly improved for all other weakly damped modes.}


 %
\begin{figure}[tb]
	\centering{
	\includegraphics[width=0.99\columnwidth]{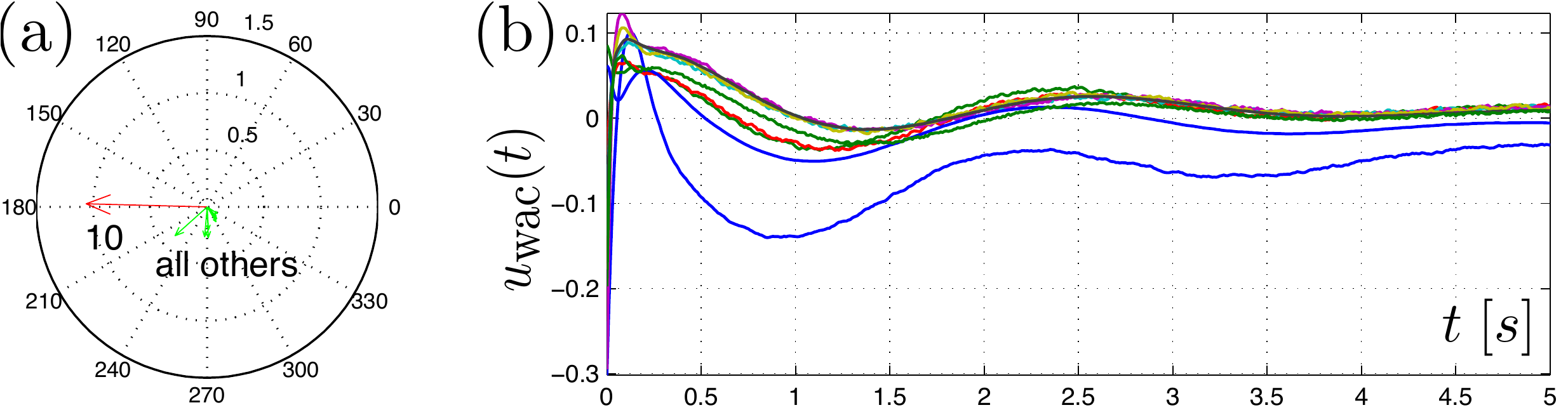}
	}
	\caption{Subfigure (a) displays the frequency components of the least damped oscillatory mode of the closed loop with local PSSs and WAC. Subfigure (b) shows the WAC signal $\subscr{u}{wac}(t)$. The initial conditions are aligned with the eigenvector of the dominant open-loop inter-area mode 1, and $\subscr{u}{wac}(t)$ is subject to additive white noise with zero mean and standard deviation~$0.01$.}
	\label{Fig: time-series NE 1}
\end{figure}

\begin{figure*}[bt]
	\centering{
	\includegraphics[width=1.99\columnwidth]{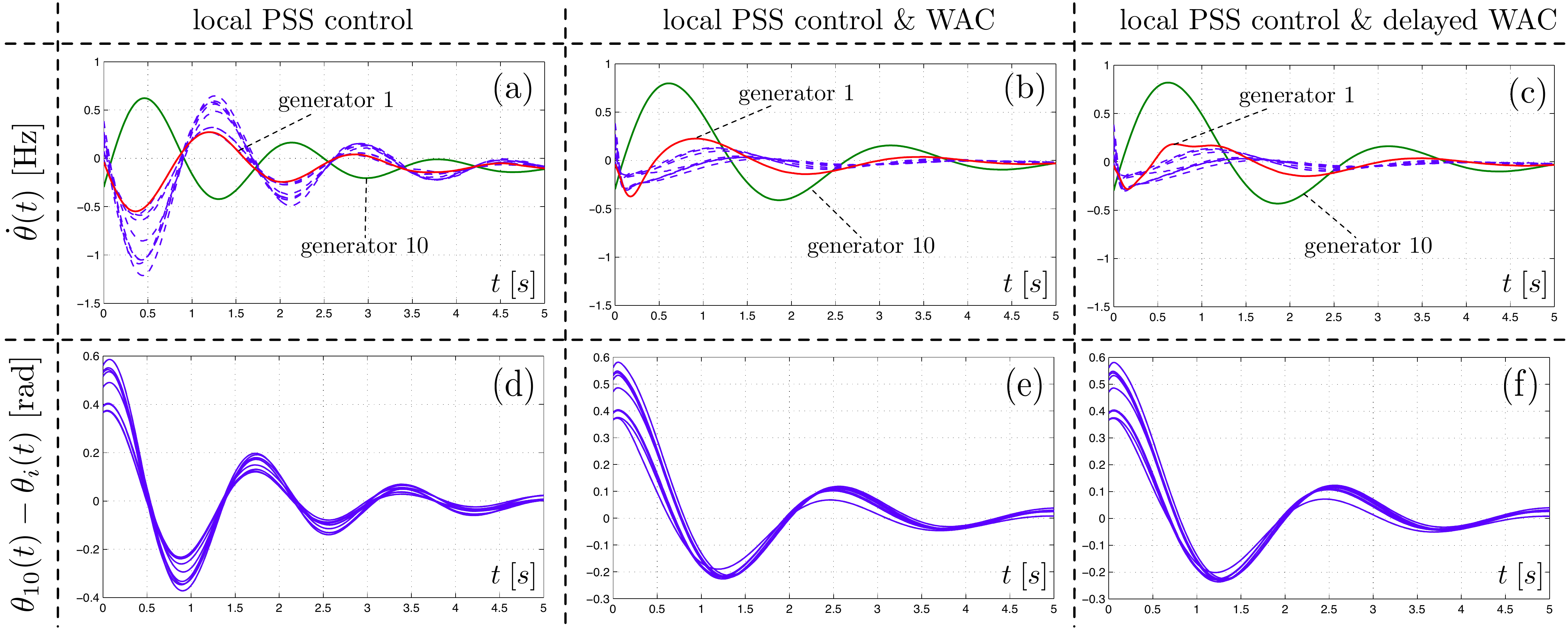}
	\caption{Time-domain simulation of the linearized model \eqref{eq: power network model} for the New England power grid: The subfigures show the generator frequencies and the angle differences (corresponding to inter-area power transfers) resulting from the use of only local PSS control, local PSS control and WAC (with the WAC signal $\subscr{u}{wac}(t)$), and local PSS control and delayed WAC (with the remote WAC signal $\subscr{u}{wac}^{\textup{(rem)}}(t)$ from \eqref{eq: decomposition of WAC} being delayed by $750$\textup{ms}).
	The initial conditions are aligned with the eigenvector of the dominant open-loop inter-area mode 1, and $\subscr{u}{wac}(t)$ is subject to additive white noise with zero mean and standard deviation~$0.01$.}
	\label{Fig: time-series NE 2}
	}
\end{figure*}

\subsection{Implementation issues, robustness, and delays}
	
As discussed in Subsection \ref{Subsection: Robustness and time delays}, multivariable phase margins are well suited robustness measures 
to account for communication delays, latencies in the SCADA network, and asynchronous measurements of generator rotor angles. Here, we investigate {how} the general multivariable phase margin of the WAC closed loop {changes with the sparsity-promoting parameter $\gamma$}; see Fig.~\ref{Fig: multivariable phase margin}. 
	\begin{figure}[bt]
	\centering{
	\includegraphics[width=0.99\columnwidth]{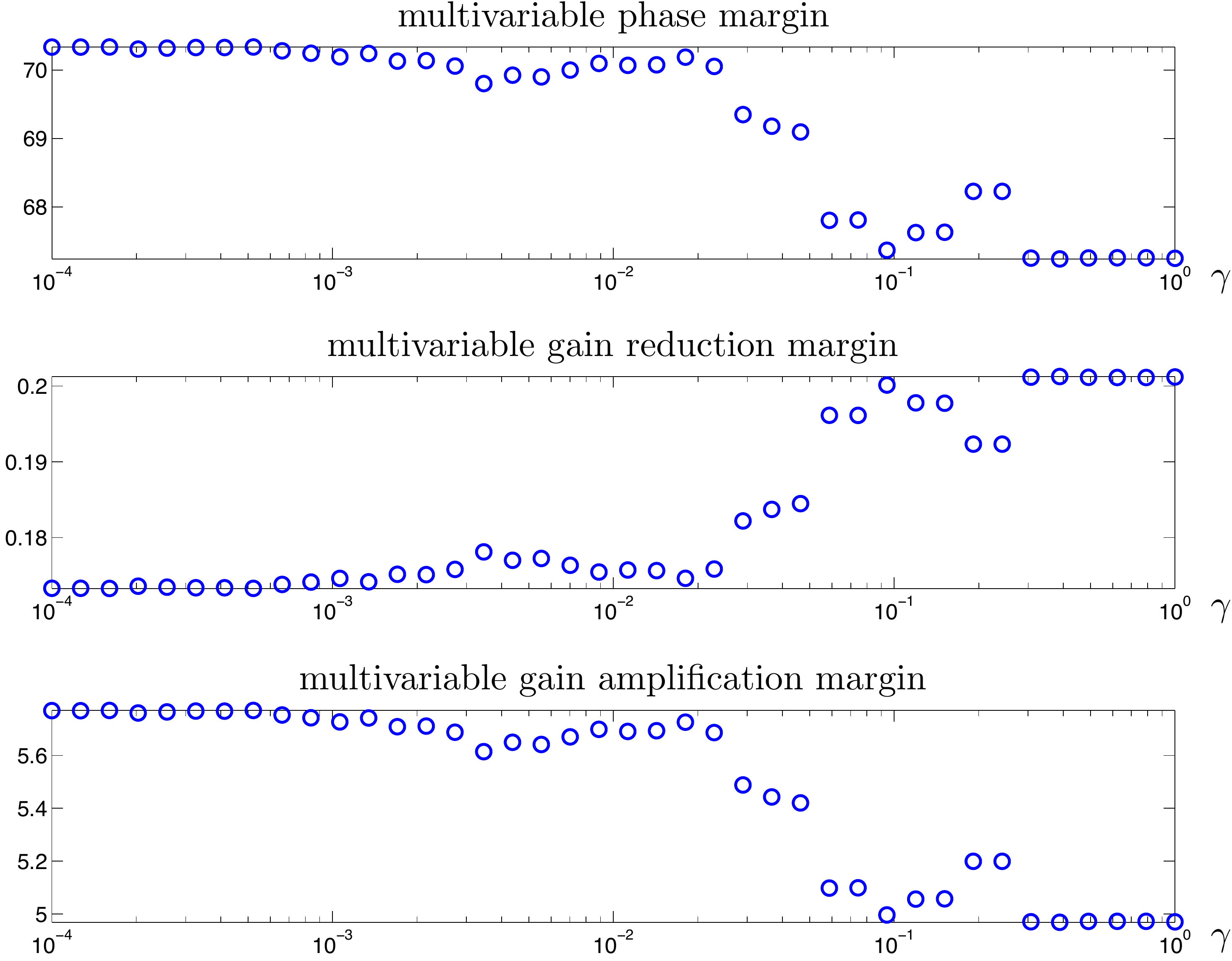}
	\caption{Multivariable phase and gain margins achieved by the {optimal} sparsity-promoting {controller} as a function of $\gamma$.}
	\label{Fig: multivariable phase margin}
	}
\end{figure}
Similar to the performance and sparsity {levels of $K^*_\gamma$} in Fig.~\ref{Fig: NE - performance}, the multivariable phase margin is nearly {a monotonic} function of $\gamma \in [10^{-4},10^{0}]$. {Over the range of examined values of $\gamma$}, the phase margin {experiences a modest drop of} about $4^{\circ}$. 
{Similar observations can be made about} the multivariable gain reduction and amplification margins displayed in Fig.~ \ref{Fig: multivariable phase margin}. They deteriorate gracefully, in {a monotonic} way, and by comparatively small amounts. Hence, each PSS gain can differ from its nominal gain (used to construct the linearized model \eqref{eq: power network model}) and independently from other PSS gains by a factor of order $[0.2,5]$ even for $\gamma=1$.
We {thus} conclude that the {controller} resulting from the sparsity-promoting optimal control problem \eqref{eq: optimal control problem} enjoys nearly the same robustness properties as the {optimal centralized controller}.

The wide-area control signal can be {further} decomposed as 
\begin{equation}
	\subscr{u}{wac}(t) 
	~ = ~
	\subscr{u}{wac}^{\textup{(loc)}}(t) 
	~ + ~
	\subscr{u}{wac}^{\textup{(rem)}}(t)\,,
	\label{eq: decomposition of WAC}
\end{equation}
where $\subscr{u}{wac}^{\textup{(loc)}}(t)$ corresponds to block-diagonal state feedback, and $\subscr{u}{wac}^{\textup{(rem)}}(t)$ corresponds to {feedback that utilizes} the remote measurement signal $\theta_{9}(t)$ {(which has to be communicated to the AVR controller at generator 1).} 
{Since the block-diagonal state feedback $\subscr{u}{wac}^{\textup{(loc)}}(t)$ depends only on generator internal states,} it can be implemented locally. The required generator state variables can be reconstructed online using appropriate filters and readily available measurements \cite{KZ-XSD-HC:13}.

For the wide-area controller~$\subscr{u}{wac}(t) = -K_{1}^{*}x(t)$,  if the local control $\subscr{u}{wac}^{\textup{(loc)}}(t)$ is absorbed in the plant, then we obtain a $58.303^{\circ}$ phase margin with respect to the remote control input $\subscr{u}{wac}^{\textup{(rem)}}(t)$. The corresponding tolerable time-delay margin is $3.964\,\textup{s}$. {From a simulation of the linearized power grid model \eqref{eq: power network model}, we observe that a delay of $750$\textup{ms} of the remote feedback signal $\subscr{u}{wac}^{\textup{(rem)}}(t)$ does not significantly affect the closed-loop performance; see Fig.~\ref{Fig: time-series NE 2}(c) and  Fig.~\ref{Fig: time-series NE 2}(f). In the delayed case, only a slight distortion of the trajectories of generator 1 (which is implementing the remote WAC signal $\subscr{u}{wac}^{\textup{(rem)}}(t)$) is noticed.} 

Additionally, we found that the information structure identified by the WAC {feedback} $\subscr{u}{wac}(t) =-K_{1}^{*}x$ is not~sensitive to the actual operating and linearization point in the dynamics \eqref{eq: non-linear DAE model}. To validate this hypothesis, we randomly altered the 
power demand at each load to create different loading conditions and operating points. Additionally, we perturbed the nominal PSS gains in \eqref{eq: PSS design}.
%
As a result, even if the PSS gains are randomly altered by factors within $[0.5,4]$, and the load power demand is randomly altered within $\pm 25\%$ of the nominal demand (leading to different linearization matrices in \eqref{eq: power network model}), 
the sparsity pattern of $K_{1}^{*}$ is identical to the one shown in Fig.~ \ref{Fig: NE - Spy(BK)}. 

We conclude that the feedback resulting from the sparsity-promoting optimal control problem \eqref{eq: optimal control problem} is not only characterized by low communication requirements and good closed-loop performance but also {by} favorable robustness margins.

\subsection{Sparsity identification and alternative control schemes}

{We note that while} the local state feedback $\subscr{u}{wac}^{\textup{(loc)}}(t)$ essentially {requires} a retuning of the local generator excitation control, the remote control signal $\subscr{u}{wac}^{\textup{(rem)}}(t)$ needs to be communicated. Recall that one of our motivations for the sparsity-promoting optimal control formulation \eqref{eq: optimal control problem} -- besides finding a stabilizing, robust, and optimal feedback -- was the {{{\em identification}} of an appropriate WAC architecture: which remote measurements need to be accessed by which controller? }

For the considered model, our sparsity-promoting {optimization framework} identified a single crucial WAC channel: the rotor angle $\theta_{9}(t)$ serves as measurement to the AVR control at generator 1. 
After having identified this WAC channel, alternative control schemes can be developed. For~instance, in absence of the local ``retuning'' control $\subscr{u}{wac}^{\textup{(loc)}}(t)$, {the WAC feedback signal} $\subscr{u}{wac}(t)$ reduces to the proportional {feedback}
 \begin{equation}
 	\subscr{u}{wac}^{\textup{(rem)}}(t) 
	~=\;
	{
	\left[ 
	\begin{array}{cccc}
	k \left(\theta_{9}(t) \, - \, \theta_{1}(t) \right) 
	&
	0 
	&
	\cdots 
	&
	0 
	\end{array}
	\right]^T.
	}
	\label{eq: remote wide-area control signal}
 \end{equation}
{Here,} $k>0$ is the $(1,9)$ element of $K_{1}^{*}$ and we {use} the local rotor angle $\theta_{1}(t)$ as {a} reference for $\theta_{9}(t)$. Hence, no absolute angle measurements are required%
 \footnote{Indeed, in the limit $\varepsilon \searrow 0$ in the cost functions \eqref{eq: state cost -- average performance}-\eqref{eq: state cost -- particular mode}, the solution to the optimal control problem \eqref{eq: optimal control problem - H2} features a rotational symmetry, and the control input requires only angular differences rather than absolute angles.},
 and the WAC control \eqref{eq: remote wide-area control signal} can be implemented, for example,  
 by integrating the frequency difference $\dot\theta_{9}(t) - \dot\theta_{1}(t)$.
 %
 The simple yet crucial proportional remote control signal \eqref{eq: remote wide-area control signal} together with an appropriate retuning of the local PSSs \eqref{eq: PSS design} then yields the WAC signal $\subscr{u}{wac}(t)$.
 
 To illustrate the utility of our {sparsity-promoting identification scheme}, we apply the proportional WAC signal \eqref{eq: remote wide-area control signal} to the full nonlinear differential-algebraic power network model \eqref{eq: non-linear DAE model} equipped with the PSSs \eqref{eq: PSS design}. To demonstrate the importance of the WAC channel $9 \to 1$, we do not re-tune the local PSSs to emulate $\subscr{u}{wac}^{\textup{(loc)}}(t)$. {We verify that the WAC signal \eqref{eq: remote wide-area control signal} does not lead to a loss of small-signal stability, and we investigate the performance of the closed-loop using time-domain simulations in Fig.~\ref{Fig: nonlinear sim}.}
 Our simulation scenario is based on a three phase fault at line $\{3,4\}$ at $0.1$\textup{s}, which is cleared at $0.2$\textup{s}. The WAC input \eqref{eq: remote wide-area control signal} is subject to a delay of {750ms}. {In analogy to the linear simulation in Fig.~\ref{Fig: time-series NE 2}, we compare the system responses with and without the WAC signal \eqref{eq: remote wide-area control signal} both with and without a 750\textup{ms} delay.} As in the linearized case, the decay rates of the frequencies are comparable but the inter-area~modes are distorted in WAC closed loop; see Fig.~\ref{Fig: nonlinear sim}(a) and \ref{Fig: nonlinear sim}(b). Due to the lack of tuned local feedback $\subscr{u}{wac}^{\textup{(loc)}}(t)$ in \eqref{eq: remote wide-area control signal}, generator 10 still slightly swings against the remaining generators in WAC closed loop. On the other hand, the implementation of the WAC signal \eqref{eq: remote wide-area control signal} by generator 1 results in a reduced  effort in providing stabilizing control signals; see Fig.~\ref{Fig: nonlinear sim}(d) and \ref{Fig: nonlinear sim}(e). {These conclusions also hold in presence of a severe delay, as seen in Fig.~\ref{Fig: nonlinear sim}(c) and \ref{Fig: nonlinear sim}(f). As compared to the linear case (including the local feedback $\subscr{u}{wac}^{\textup{(loc)}}(t)$) shown in Fig.~\ref{Fig: time-series NE 2}, the nonlinear performance (without local feedback $\subscr{u}{wac}^{\textup{(loc)}}(t)$) is not insensitive to delays and slightly degrades. To further account for communication delays, a local control $\subscr{u}{wac}^{\textup{(loc)}}(t)$ can beed added, and the remote proportional WAC feedback signal \eqref{eq: remote wide-area control signal} can be fed through a delay-compensating filter as in \cite{PZ-DYY-KWC-GWC:12}.}
 
 \begin{figure*}[t]
	\centering{
	\includegraphics[width=1.99\columnwidth]{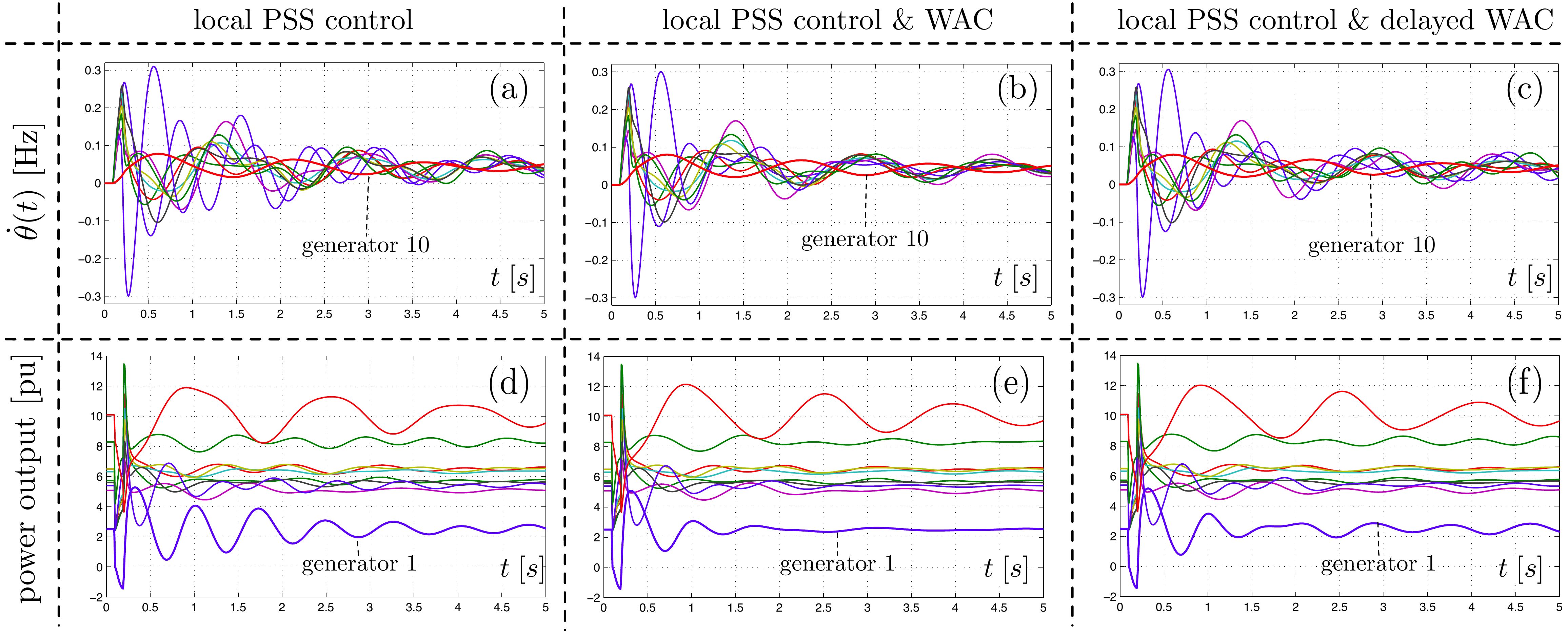}
	\caption{{Time-domain simulation of the full nonlinear differential-algebraic model \eqref{eq: non-linear DAE model} for the New England power grid:  The subfigures show the generator frequencies and active power outputs using only local PSS control and using the remote WAC signal $\subscr{u}{wac}^{\textup{(rem)}}(t)$ both without delay and with a 750$\textup{ms}$~delay. Initially, the system is at steady state, and a three phase fault at line $\{3,4\}$ at $0.1$\textup{s}, which is cleared at $0.2$\textup{s}.}}
	\label{Fig: nonlinear sim}
	}
\end{figure*}

We conclude that the simple proportional {feedback} \eqref{eq: remote wide-area control signal} identified by the sparsity-promoting {framework} significantly improves the closed-loop performance. {We note that the  proportional WAC strategy \eqref{eq: remote wide-area control signal} merely serves as {an illustration of the ability of our method to identify} a suitable WAC architecture. Now that the crucial WAC channel $9 \to 1$ has been identified, more sophisticated single-input-single-output control strategies can be developed via alternative means.} 

\section{Conclusions}
\label{Section: Conclusions}

We proposed a novel approach to {wide-area control} of inter-area oscillations. We followed a recently {introduced paradigm} to sparsity-promoting optimal control \cite{MF-FL-MRJ:13-updated}. Our performance objectives were inspired by the well-known slow coherency theory. We {validated} the proposed control strategy {on the  IEEE 39 New England power grid model}, which demonstrated nearly optimal performance, low communication requirements, and robustness of {the sparsity-promoting controller}.

{Of course, further scrutiny of the proposed sparsity-promoting framework on other power system models is needed. We note that both the structure and the performance achieved by our sparsity-promoting controller are problem-dependent: they are influenced by the model of the underlying power system as well as by the selection of state and control weights. For the considered IEEE 39 New England model, we have illustrated the insensitivity of the identified controller structure on the operating point and the graceful performance degradation relative to the optimal centralized controller. In our future effort, we will examine if these properties also hold for other power systems and investigate how the resulting controller structure is influenced by the systemÕs dynamics. }


Our initial results are very promising, and we are currently working on further extensions. The authors also envision an application of the proposed control strategy to enhance the performance of purely decentralized primary and secondary control schemes \cite{DK-SB-MC:12} with supplementary WAC. 
Since sparsity is a function of the state space coordinates, the identified control scheme may not be sparse in other coordinates. Hence, an interesting question is to identify state space representations of power system models, which are amenable to a sparsity-promoting control design. For example, it is of interest to extend the proposed design methodology to sparse, differential-algebraic, and structure-preserving models. 
{Finally, the authors also aim to gain insight into fundamental performance limitations induced by localized and wide-area control strategies.}

\section*{Acknowledgments}

The authors would like to thank Sairaj Dhople for helpful comments on an early version of this paper and Scott Ghiocel for helpful discussions regarding the Power Systems Toolbox.


\begin{appendix}

\subsection{Algorithmic methods for sparsity-promoting optimization}
\label{Subsection: Algorithmic Methods}

We briefly summarize the algorithmic approach to the optimization problem  \eqref{eq: optimal control problem - H2} and refer to \cite{MF-FL-MRJ:13-updated} for further details:
\begin{enumerate}
	
	\item[(i)] {\bf Warm-start and homotopy:} The optimal control problem~\eqref{eq: optimal control problem - H2} is solved by tracing a homotopy path that starts at the optimal centralized controller with $\gamma = 0$ and continuously increases $\gamma$ until the desired value $\gamma_{\mathrm{des}}$; 
	
	\item[(ii)] {\bf ADMM:} For each value of $\gamma \in [0, \, \gamma_{\mathrm{des}}]$, the optimization problem \eqref{eq: optimal control problem - H2} is solved iteratively using ADMM; 
	
	\item[(iii)] {\bf Updates of weights:} In each step of ADMM, the weights are updated as $w_{ij}  = 1/(|K_{ij}| + \varepsilon)$ with $\varepsilon>0$.
	We have conducted $5$ update steps with $\varepsilon = 10^{-3}$; and
	
	\item[(iv)] {\bf Polishing:} Once the desired sparsity pattern $\mc K$ is identified, a structured optimal control problem is solved:%
	\end{enumerate}
	\begin{align*}
		\begin{split}
		&\minimize \quad 	J_{\gamma}(K) =\;  \textup{trace}\bigl(B_{1}^{T}  P B_{1}\bigr)  \,,
		\\
	    &\textup{subject to} \;\quad K \in \mc K \,, 
	    \\& \quad  
	    \bigl(A-B_{2}K)^{T} P + P(A-B_{2}K) = - (Q+K^{T}RK) \,.
		\end{split}
	\end{align*}
This iterative approach is convergent under a local convexity assumption, and stability of $A-B_{2}K$ is guaranteed.
%
The algorithms developed in \cite{MF-FL-MRJ:13-updated} have been implemented in {\em MATLAB}. The associated software {as well as the numerical power network data} used in the present paper can be downloaded at 
{\small \tt www.ece.umn.edu/users/mihailo/software/lqrsp/}.
\end{appendix}



\bibliographystyle{IEEEtran}
\bibliography{alias,Main,FB,New} 


\end{document}